\documentclass[12pt,leqno]{article}

\usepackage{amsmath,amssymb}
\def\R{\Bbb R}
\def\C{\Bbb C}

\def\Z{\Bbb Z}
\def\T{\Bbb T}

\def\O{\Omega}

\def\supp{{\rm supp}}

\def\cL{{\cal L}}
\def\cC{{\cal C}}

\def\cT{{\cal T}}
\def\cI{{\cal I}}
\def\cP{{\cal P}}

\def\Linfg{L^\infty (G)}

\def\spec{{\rm spec}}
\def\cMT{{\cal M}_T(\Sigma)}

\begin{document}
\title{Decomposition of analytic measures on groups and measure spaces}
\author{Nakhl\'e Asmar and Stephen Montgomery-Smith\\
Department of Mathematics\\
University of Missouri\\
Columbia, MO 65211\\
{\tt nakhle@math.missouri.edu}\\
{\tt stephen@math.missouri.edu}}
\date{}

\maketitle

\begin{center}
{\em Dedicated to the memory of Edwin Hewitt}
\end{center}

\newtheorem{defin}{Definition}[section]
\newtheorem{thm}[defin]{Theorem}
\newtheorem{ex}[defin]{Example}
\newtheorem{prop}[defin]{Proposition}
\newtheorem{lemma}[defin]{Lemma}
\newtheorem{schol}[defin]{Scholium}
\newtheorem{rem}[defin]{Remarks}
\newtheorem{cor}[defin]{Corollary}

\begin{abstract}
In this paper, we consider an arbitrary locally compact
abelian group $G$, with an ordered dual group $\Gamma$,
acting on a space of measures.  
Under suitable conditions, we define the notion of 
analytic measures using the representation of $G$ and the order on $\Gamma$.  
Our goal is to study analytic measures by applying a new  transference principle for subspaces of measures, along with  results from probability and Littlewood-Paley theory.  As a consequence,  
we will derive new properties of analytic measures 
 as well as extensions of previous work of Helson and Lowdenslager, de Leeuw and Glicksberg, and Forelli.

A.M.S. Subject Classification:  43A17, 43A32,  

Keywords: orders, transference, measure space, sup path attaining, 
F.\&M.\ Riesz Theorem
\end{abstract}

\section{Introduction}

This paper is essentially providing a new approach to generalizations of 
the F.\&M.~Riesz Theorems, for example, such results as that  
of Helson and Lowdenslager \cite{hl1, hl2}.  They showed that if 
$G$ is a compact abelian group with ordered dual, and if $\mu$ is an
{\em analytic\/} measure (that is, its Fourier transform is supported
on the positive elements of the dual), then it follows that the singular
and absolutely continuous parts (with respect to the Haar measure) are also
analytic.

Another direction is that provided by Forelli \cite{forelli}
(itself a generalization of the result of 
de Leeuw and Glicksberg \cite{deleeuwglicksberg}), 
where one has an
action of the real numbers $\R$ acting on a locally compact topological
space $\Omega$, and a Baire measure $\mu$ on $\Omega$ that is {\em analytic\/}
(in a sense that we make precise below) with respect to the action.  Then
again, the singular and absolutely continuous parts of $\mu$ (with respect
to any so called quasi-invariant measure) are also analytic.

Indeed common generalizations of both these ideas have been provided, for
example, by Yamaguchi \cite{yama}, considering the action of any locally
compact abelian group with ordered dual, on a locally compact topological
space.  For more generalizations we refer the reader to
 Hewitt, Koshi, and Takahashi \cite{hkt}.

In the paper \cite{amss}, a new approach to proving these kinds of results
was given, providing a transference principle for spaces of measures.  In
that paper, the action was from a locally compact abelian group into a
space of isomorphisms on the space of measures of a sigma algebra.  
A primary requirement that the action had to satisfy was what was called
{\em sup path attaining}, a property that was satisfied, for example, by
the setting of Forelli (Baire measures on a locally compact topological 
space).  Using this transference principle, the authors were able to 
give an extension and a new proof of Forelli's result.  This was obtained
by using a Littlewood-Paley decomposition of an analytic measure.

In this paper we wish to continue this process, applying this same 
transference principle to provide the common generalizations of the
results of Forelli and Helson and Lowdenslager.  What we provide in
this paper is essentially a decomposition of an analytic measure as 
a sum of martingale differences with respect to
a filtration defined by the order.  For each martingale difference, the
action of the group can be described precisely by a certain action of the
group of real numbers, and so we can appeal to the results of \cite{amss}.

In this way, we can reach the following generalization
(see Theorem~\ref{application1} below): if $\cal P$ is any
bounded operator on the space of measures that commutes with the action
(as does, for example, taking the singular part), and if $\mu$ is an
analytic measure, then ${\cal P}\mu$ is also an analytic measure.

In the remainder of the introduction, we will establish our notation,
including the notion of sup path attaining,
and recall the transference principle from \cite{amss}.  In Section~2,
we will describe orders on locally compact abelian groups, including
the extension of Hahn's Embedding Theorem provided in \cite{ams}.
In Section~3, we define the notions of analyticity.  This somewhat
technical section continues into Section~4, which examines the role
of homomorphism with respect to analyticity.  The technical results
basically provide proofs of what is believable, and so may be skipped
on first reading.  It will be seen that the concept of sup path attaining
comes up again and again, and may be seen to be an integral part of all
our proofs.

In Section~5, we are ready to present the decomposition of
analytic measures.  This depends heavily on transference of martingale 
inequalities of Burkholder and Garling, and then using the fact 
 that weakly unconditionally summing series are
unconditionally summing in norm for any series in a space of measures \cite{bp}.
In Section~6, we then give applications of this decomposition,
giving the generalizations that we alluded to above.

Throughout $G$ will denote a 
locally compact abelian group with dual group
$\Gamma$.  
The symbols $\Z$, $\R$ and $\C$ 
denote the integers, the real and complex numbers, respectively.  
If $A$ is a set, we denote the indicator
function of $A$ by $1_A$.
For $1\leq p<\infty$,
the space of Haar measurable functions $f$ on $G$ with
$\int_G|f|^p dx<\infty$ will be denoted by 
$L^p(G)$.  The space of essentially
bounded functions on $G$ will be denoted by 
$L^\infty(G)$.  The expressions ``locally null''
and ``locally almost everywhere'' will have the same meanings as
in \cite[Definition (11.26)]{hr1}.

Let $\cC_0(G)$ denote
the Banach space of continuous functions on
$G$ vanishing at infinity. 
The space of all complex regular Borel measures 
on $G$,
denoted by $M(G)$, consists of all complex measures  
arising from bounded linear functionals on $\cC_0(G)$.

Let $(\O, \Sigma)$ denote a
measurable space, 
where $\O$ is a set and $\Sigma$ is a
sigma algebra of subsets of $\O$.  Let $M(\Sigma)$ denote the 
Banach space of complex measures on $\Sigma$ with the
total variation norm, and let
$\cL^\infty(\Sigma)$ denote the space of measurable
bounded functions on $\Omega$.  

Let
$T:\ t\mapsto T_t$ denote a representation of $G$
by isomorphisms of $M(\Sigma)$.  
We suppose that $T$ is uniformly bounded,
i.e., there is a positive constant $c$ such that
for all $t\in G$, we have
\begin{equation}
\|T_t\|\leq c .
\label{uniformlybded}
\end{equation}
\begin{defin}
A measure $\mu\in M(\Sigma)$ is 
called weakly measurable (in symbols, $\mu\in{\cal M}_T(\Sigma)$)
if for every $A\in \Sigma$ 
the mapping $t\mapsto T_t\mu(A)$ is 
Borel measurable on $G$.  
\label{weakmble}
\end{defin}

Given a measure $\mu\in \cMT$ and a Borel measure 
$\nu \in M(G)$, we define the `convolution' 
$\nu*_T\mu$ on $\Sigma$ by
\begin{equation}
\nu*_T\mu (A)=\int_G T_{-t}\mu(A) d\nu(t)
\label{Tconv}
\end{equation}
for all $A\in\Sigma$.

We will assume throughout this paper
that the representation $T$ commutes with the 
convolution (\ref{Tconv}) in the following sense:
 for each $t\in G$, 
\begin{equation} 
 T_t(\nu*_T\mu)=\nu*_T(T_t\mu).
\label{commut}
\end{equation}
Condition (\ref{commut}) holds if, for example,
for all $t\in G$, the adjoint of $T_t$ maps 
$\cL^\infty(\Sigma)$ into itself.  In symbols,
\begin{equation}
T_t^*: \cL^\infty(\Sigma) \rightarrow \cL^\infty(\Sigma).
\label{adjointT}
\end{equation}
For proofs we refer the reader to 
\cite{ams2}.  Using (\ref{uniformlybded}) and 
(\ref{commut}), it can be shown that $\nu*_T\mu$ is a measure in $\cMT$,
\begin{equation}
\|\nu*_T\mu\|\leq  c\|\nu\|\|\mu\|,
\label{normofconv}
\end{equation}
where $c$ is as in 
(\ref{uniformlybded}), and 
\begin{equation} 
\sigma*_T(\nu*_T\mu)=(\sigma*\nu)*_T\mu,
\label{associative}
\end{equation}
for all $\sigma , \nu \in M(G)$ and $\mu \in \cMT$
(see \cite{ams2}).

\begin{defin}
A representation $T=(T_t)_{t\in G}$ of a locally compact abelian group
$G$ in $M(\Sigma)$ is said to be sup path attaining
if it is uniformly bounded, satisfies property (\ref{commut}), and if there 
is a constant $C$ such that for every weakly 
measurable $\mu\in {\cal M}_T(\Sigma)$ we have
\begin{equation}
\| \mu\| \leq C\sup \left\{
{\rm ess\ sup}_{t\in G} \left|
\int_\O h d (T_t\mu)
\right|
:\ \ h\in \cL^\infty(\Sigma),\ \|h\|_\infty\leq 1
\right\}.
\label{ineqhypa}
\end{equation}
\label{def hypa}
\end{defin}

The fact that the mapping $t\mapsto \int_\O h d (T_t\mu)$ is 
measurable is a simple consequence of
the measurability of the mapping $t\mapsto T_t\mu(A)$ for every
$A\in\Sigma$.  

In \cite{amss} were provided many examples of sup path attaining
representations.  Rather than give this same list again, we give
a couple of examples of particular interest.

\begin{ex}
{\rm (a) (This is the setting of Forelli's Theorem.)
Let $G$ be a locally compact abelian group, and
$\Omega$ be a locally compact topological
space.  Suppose that $\left( T_t\right)_{t\in G}$ is a
group of homeomorphisms of $\Omega$ onto itself
such that the mapping
$$(t,\omega)\mapsto T_t\omega$$
is jointly continuous.
Then the space of Baire measures on $\Omega$, that is, the minimal
sigma algebra such that compactly supported continuous functions are
measurable, is sup path attaining under the action
$T_t\mu(A)=\mu(T_t(A))$, where $T_t(A)=\{T_t\omega:\ \omega\in A\}$.
(Note that all Baire measures are weakly measurable.)

\noindent
(b)
Suppose that $G_1$ and $G_2$ are locally compact abelian groups and that
 $\phi:\ G_2\rightarrow G_1$ is a continuous homomorphism.
Define an action of $G_2$ on $M(G_1)$ (the regular Borel measures
on $G_1$) by translation by $\phi$.
Hence, for $x\in G_2, \mu\in M(G_1)$, and any
Borel subset $A\subset G_1$, let $T_x\mu(A)=\mu(A+\phi(x))$.
Then every $\mu\in M(G_1)$ is weakly measurable, and the
representation
is sup path attaining with constants $c = 1$ and $C = 1$.
}
\label{exhypa}
\end{ex}

\begin{prop}
Suppose that $T$ is sup path attaining 
and $\mu$ is weakly measurable such that for 
every $A\in \Sigma$ we have
$$T_t\mu(A)=0$$
for locally almost all $t\in G$.  Then $\mu=0$.
\label{prop hypa}
\end{prop}
The proof is immediate (see \cite{ams2}).\\

We now 
recall some basic definitions from spectral theory.

 
If $I$ is an ideal in $L^1(G)$, let
$$Z(I)=\bigcap_{f\in I}
\left\{
\chi\in\Gamma:\ \ \widehat{f}(\chi)=0
\right\}.$$
The set $Z(I)$ is called the zero set of $I$.
For a weakly measurable $\mu\in M(\Sigma)$, let
$$\cI (\mu)=\{f\in L^1(G):\ \ f*_T\mu =0\}.$$
When we need to be specific about the 
representation, we will use 
the symbol $\cI_T (\mu)$ instead of $\cI (\mu)$.

Using properties of the convolution $*_T$,
it is straightforward to show that $\cI(\mu)$ is a closed ideal
in $L^1(G)$.
\begin{defin}
The $T$-spectrum of a 
weakly measurable $\mu\in \cMT$ is defined by 
\begin{equation}
{\rm spec}_T (\mu)= \bigcap_{f\in \cI(\mu)}
\left\{
\chi\in\Gamma:\ \ \widehat{f}(\chi)=0
\right\}=Z(\cI(\mu)).
\label{specsbt}
\end{equation}
\label{Tspectrum}
\end{defin}
If $S\subset \Gamma$, let
$$L_S^1=L_S^1(G)=\left\{f\in L^1(G):\ \widehat{f}=0\ \mbox{outside of}\ S\right\}\,.$$
In order to state the main transference result, we introduce one more definition.

\begin{defin}
A subset $S\subset\Gamma$ is a $\cT$-set if, given any compact 
$K\subset S$, each neighborhood of $0\in\Gamma$ contains a nonempty open set
$W$ such that $W+K\subset S$. 
\label{t-set}
\end{defin}

\begin{ex}\label{s-set}
{\rm (a)  If $\Gamma$ is a locally compact abelian group, then any open subset of $\Gamma$ is a $\cT$-set.  In particular, if $\Gamma$ is discrete then every subset of $\Gamma$ is a $\cT$-set.\\
(b)  The set $\left[ a,\infty\right. )$ is a $\cT$-subset of $\R$, for all $a\in\R$.  \\
(c)  Let $a\in\R$ and $\psi:\ \Gamma \rightarrow \R$ be a continuous homomorphism.  Then 
$S=\psi^{-1}([a,\infty))$ is a $\cT$-set.\\
(d)  Let $\Gamma=\R^2$ and $S=\{(x,y):\ y^2\leq x\}$.  
Then $S$ is a $\cT$-subset of $\R^2$ such that there is no nonempty open set $W\subset \R^2$ such that $W+S\subset S$. 
}
\end{ex}

The main result of \cite{amss} is the following transference
theorem.
\begin{thm}
Let $T$ be a sup path attaining representation of a locally compact abelian group
$G$ by isomorphisms of $M(\Sigma)$ and let $S$ be a $\cT$-subset of
$\Gamma$.
Suppose that $\nu$ is a measure in $M(G)$ such that
\begin{equation}
\|\nu*f\|_1\leq \|f\|_1
\label{hyptransference2}
\end{equation}
for all $f$ in $L_S^1(G)$.  Then for every 
weakly measurable $\mu \in M(\Sigma)$ with
${\rm spec}_T ( \mu )\subset S$
we have
\begin{equation}
\|\nu*_T\mu\|\leq c^3 C \|\mu\|,
\end{equation}
where $c$ is as in (\ref{uniformlybded}) and 
$C$ is as in (\ref{ineqhypa}).
\label{trans-thm}
\end{thm}

\section{Orders on locally compact abelian groups}

An order $P$ on $\Gamma$ is a subset
that satisfies the three axioms:
$P+P\subset P$; $P\cup (-P)=\Gamma$; and $P\cap (-P)-\{0\}$.
We recall from \cite{ams} the following property of orders.

\begin{thm}
Let $P$ be a measurable order on $\Gamma$.  There 
are a totally ordered set $\Pi$ with largest element 
$\alpha_0$; a chain of subgroups $\{C_\alpha\}_{\alpha\in\Pi}$
of $\Gamma$; 
and a collection of continuous real-valued homomorphisms 
$\{\psi_\alpha\}_{\alpha\in\Pi}$ on $\Gamma$ such that:\\
(i)\ \ for each $\alpha\in\Pi$, $C_\alpha$ is an 
open subgroup of $\Gamma$;\\
(ii)\ \ $C_\alpha\subset C_\beta$ if $\alpha > \beta$.\\ 
Let $D_\alpha=\{\chi\in C_\alpha:\ \psi_\alpha(\chi)=0\}$.  Then,
 for every $\alpha\in \Pi$,\\
(iii) $\psi_\alpha(\chi)>0$ for every 
	$\chi\in P\cap (C_\alpha\setminus D_\alpha)$,\\
(iv)    $\psi_\alpha(\chi)<0$ for every $\chi\in 
	(-P)\cap (C_\alpha\setminus D_\alpha).$\\
(v)\ \ When $\Gamma$ is discrete, $C_{\alpha_0}=\{0\}$; and 
when $\Gamma$ is not discrete, $D_{\alpha_0}$ 
has empty interior and is locally null.
\label{structureorder}
\end{thm}

When $\Gamma$ is discrete, Theorem \ref{structureorder} 
can be deduced from
the proof of Hahn's Embedding Theorem for
orders (see \cite[Theorem 16, p.59]{fu}). 
The general case treated in Theorem \ref{structureorder}
accounts for the measure theoretic aspect of orders.
The proof is
based on the study of orders of Hewitt and Koshi \cite{hk}.

For $\alpha\in \Pi$\ with $\alpha\neq\alpha_0$, let
\begin{eqnarray}
S_\alpha\equiv P\cap (C_\alpha \setminus D_\alpha)
&=& \left\{
\chi\in C_\alpha\setminus D_\alpha:\ \ \psi_\alpha(\chi)\geq 0\right\}
\label{alphaslice}\\
&=&\left\{
\chi\in C_\alpha :\ \ \psi_\alpha(\chi)> 0\right\}.
\label{alphaslice2}
\end{eqnarray}
For $\alpha=\alpha_0$, set
\begin{equation}
S_{\alpha_0}=\left\{
\chi\in C_{\alpha_0}:\ \ \psi_{\alpha_0}(\chi)\geq 0\right\}.
\label{alpha0slice}
\end{equation}
Note that when $\Gamma$ is discrete, $C_{\alpha_0}=\{0\}$, and so $S_{\alpha_0}=\{0\}$ in this case. 

If $A$ is
a subset of a topological space, we will use $\overline{A}$ and
$A^\circ$ to denote the closure, respectively, the interior
of $A$. 

\begin{rem}
{\rm  
(a)  It is a classical fact that a group $\Gamma$ can be ordered if
and only if it is torsion-free.
Also, an order on $\Gamma$ is 
any maximal positively linearly independent set.  
Thus, orders abound in torsion-free abelian groups, as
they can be constructed using Zorn's Lemma to obtain a maximal positively
linearly independent set.  
(See \cite[Section 2]{hk}.)  
However, if we ask for 
measurable orders, then we are 
restricted in many ways in the choices 
of $P$ and also the topology on $\Gamma$.  As shown in \cite{hk}, 
any measurable order on $\Gamma$ has nonempty interior.
Thus, for example, while there are infinitely many orders
on $\R$, only two are Lebesgue measurable:  
$P=[0,\infty [$, and
$P=]-\infty,0]$.  
It is also shown in \cite[Theorem (3.2)]{hk} that any
order on an infinite compact torsion-free abelian group
is non-Haar measurable.  This effectively shows
that if $\Gamma$ contains a Haar-measurable order $P$,
and we use the structure theorem for locally compact abelian
groups to write $\Gamma$ as $\R^a\times \Delta$, where 
$\Delta$ contains a compact open subgroup
\cite[Theorem (24.30)]{hr1}, then
either $a$ is a positive integer, or $\Gamma$ is
discrete.  (See \cite{ams}.)  \\
(b)  The subgroups $(C_\alpha)$ are characterized as being the principal
convex subgroups in $\Gamma$ and for each $\alpha\in \Pi$,
we have 
$$D_\alpha=\bigcup_{\beta>\alpha}C_\beta.$$
Consequently, we have $C_\alpha\subset D_\beta$ if 
$\beta <\alpha$.   
By construction, the sets $C_\alpha$ are open.  For $\alpha< \alpha_0$, the subgroup $D_\alpha$ has
 nonempty interior, since it contains $C_\beta$, with
$\alpha<\beta$.  Hence for $\alpha\neq \alpha_0$,
$D_\alpha$ is open and closed.  Consequently, for 
$\alpha\neq \alpha_0$,
$C_\alpha\setminus D_\alpha$ is open and closed.

(c)  Let $\psi:\ \Gamma_1\rightarrow \Gamma_2$ be a continuous 
homomorphism between two ordered groups.  We say that
$\psi$ is order-preserving if $\psi(P_1)\subset P_2$.  
Consequently, if 
$\psi$ is continuous and order preserving, then
$\psi(\overline{P_1})\subset \overline{P_2}$.

For each $\alpha\in \Pi$, let $\pi_\alpha$ denote the
quotient homomorphism $\Gamma\rightarrow \Gamma/C_\alpha$. 
Because $C_\alpha$ is a principal subgroup, we can define an order
on $ \Gamma/C_\alpha$ by setting $\psi_\alpha(\chi)\geq 0\Longleftrightarrow
\chi\geq 0$.  Moreover, the principal convex subgroups
in $\Gamma/C_\alpha$ are precisely the images by 
$\pi_\alpha$ of the principal convex subgroups of 
$\Gamma$ containing $C_\alpha$.  (See \cite[Section 2]{ams}.)

}
\label{remarkstructureorder}
\end{rem}

We end this section with a useful property of orders.

\begin{prop}
Let $P$ be a measurable order on $\Gamma$.
Then
$\overline{P}$ is a 
${\cal T}$-set.
\label{ptset}
\end{prop}

\noindent
{\bf Proof.}\quad If $\Gamma$ is discrete, there is nothing to prove.  If $\Gamma$ is not discrete, the subgroup
$C_{\alpha_0}$ is open and nonempty.  Hence the set
$  C_{\alpha_0}\cap \{\chi\in \Gamma:\psi_{\alpha_0}
(\chi)>0\}$ is nonempty, with $0$ as a limit
point.  Given an open nonempty neighborhood $U$ of $0$,
let 
$$W=U\cap C_{\alpha_0}\cap \{\chi\in \Gamma:\psi_{\alpha_0}
(\chi)>0\}.$$ 
Then $W$ is a nonempty subset of $U\cap P$.
Moreover, it is easy to see that 
$W+\overline{P}\subset P\subset \overline{P}$,
and hence 
$\overline{P}$ is a 
${\cal T}$-set.

\section{Analyticity}

We continue with the notation of the previous section.
Using the order structure on $\Gamma$ we define
some classes of analytic functions on $G$:
\begin{eqnarray}
H^1(G)&=&\left\{
f\in L^1(G):  \widehat{f}=0\ {\rm on}\ (-P)\setminus\{0\}
\right\};
\label{h1g}\\
H^1_0(G)&=&\left\{
f\in L^1(G):  \widehat{f}=0\ {\rm on}\ -P
\right\};\label{h10g}
\end{eqnarray}
and
\begin{equation}
H^\infty(G)=\left\{
f\in L^\infty(G):  \int_G f(x)g(x)dx=0 \ {\rm for\ all}\ g\in H^1_0(G)
\right\}.
\label{hinfg}
\end{equation}
We clearly have
$$
H^1(G)=\left\{
f\in L^1(G):  \widehat{f}=0\ {\rm on}\ \overline{(-P)\setminus\{0\}}
\right\}.$$
We can now give the definition of analytic measures in 
$\cMT$.

\begin{defin}
Let $T$ be a sup path attaining representation
of $G$ by isomorphisms of $M(\Sigma)$.  
A measure $\mu\in \cMT$ is 
called weakly analytic if the mapping $t\mapsto T_t\mu(A)$ 
is in $H^\infty(G)$ for every $A\in\Sigma$.
\end{defin}

\begin{defin}
Recall the $T$-spectrum of a 
weakly measurable $\mu\in \cMT$, 
\begin{equation}
{\rm spec}_T (\mu)= \bigcap_{f\in \cI(\mu)}
\left\{
\chi\in\Gamma:\ \ \widehat{f}(\chi)=0
\right\}.
\label{spect}
\end{equation}
A measure $\mu$ in $\cMT$ is called $T$-analytic if
${\rm spec}_T (\mu)\subset \overline{P}$.

\label{weakly-analytic}
\end{defin}

That the two definitions of analyticity are
equivalent will be shown later in this section.

Since $\cI(\mu)$ is translation-invariant,
it follows readily that 
for all $t\in G$,
$$\cI(T_t\mu)=\cI(\mu),$$
and hence
\begin{equation}
{\rm spec}_T (T_t(\mu))={\rm spec}_T (\mu).
\label{feb.4.95}
\end{equation}


We now recall several 
basic results from 
spectral theory of bounded functions that will be needed in the sequel.  Our reference is \cite[Section 40]{hr2}.  
If 
$\phi$ is in $\Linfg$, write $\left[ \phi\right]$ 
for the smallest weak-* closed translation-invariant subspace of $\Linfg$ 
containing $\phi$, and let $\cI([\phi])=\cI (\phi)$ 
denote the closed translation-invariant ideal in $L^1(G)$:\\
$$ \cI (\phi)=\{f\in L^1(G): f*\phi=0\}.$$
It is clear that 
$ \cI (\phi)=\{f\in L^1(G): f*g=0, \forall g\in \left[\phi\right]\}$. 
The spectrum of $\phi$, denoted by 
$\sigma \left[\phi\right]$, 
is the set of all continuous characters of $G$ 
that belong to $\left[\phi\right]$.  
This closed subset of $\Gamma$ is also given by 
\begin{equation}
\sigma \left[\phi\right]=Z(\cI(\phi)). 
\label{spec}
\end{equation}
(See \cite[Theorem (40.5)]{hr2}.)  

Recall that a closed subset $E$ of $\Gamma$ is 
a set of spectral synthesis for $L^1(G)$, or an $S$-set, if and only if
$\cI([E])$ is the only ideal in $L^1(G)$ whose zero set
is $E$.  
 
There are various equivalent definitions of $S$-sets.
Here is one that we will use at several occasions.\\
{\em A set 
$E\subset \Gamma$ is an $S$-set if and only if
every essentially bounded function $g$ in $L^\infty(G)$
with $\sigma[g]\subset E$ is the weak-* limit 
of linear combinations of characters from $E$.}\\  
(See \cite[(40.23) (a)]{hr2}.)  This has 
the following immediate consequence.

\begin{prop}
Suppose that $B$ is an $S$-set, $g\in L^\infty(G)$,
 and 
$\spec (g)\subset B$.  
(i)  If $f$ is in $L^1(G)$ and 
$\widehat{f}=0$ on $B$, then
$f*g(x)=0$ for all $x$ in $G$.  In particular, 
$$\int_G f(x) g(-x)\, d\, x =0.$$
(ii)  If $\mu$ is a measure in $M(G)$
with $\widehat{\mu}=0$ on $B$, then
$\mu*g(x)=0$ for almost all $x$ in $G$.
\label{annihilate}
\end{prop}

{\bf Proof.}\quad Part (i) is a simple consequence of
\cite[Theorems (40.8) and (40.10)]{hr2}.
We give a proof for the sake of completeness.
Write $g$ as the weak-* limit of 
trigonometric polynomials,
$\sum_{\chi\in E} a_\chi \chi(x)$,
 with characters in $E$.  Then
\begin{eqnarray*}
\int_G f(x)g(y-x)\, d\,x &=& \lim \int_G 
\sum_{\chi\in E} a_\chi \chi(y)f(x)\chi(-x)\, d\,x\\
&=&\lim  
\sum_{\chi\in E} a_\chi \chi(y) \widehat{f}(\chi)=0
\end{eqnarray*} 
since $\widehat{f}$ vanishes on $E$.

To prove (ii), assume that $\mu*g$ is not 0 a.e..  
Then, there is $f$ in $L^1(G)$ such that
$f*(\mu*g)$ is not 0 a.e..  But this contradicts (i),
since $f*(\mu*g)=(f*\mu)*g$, $f*\mu$ is in $L^1(G)$,
and $\widehat{f*\mu}=0$ on $B$.

The following is a converse of sorts of Proposition
\ref{annihilate} and follows easily from definitions.

\begin{prop}
Let $B$ be a nonvoid closed subset of $\Gamma$.
Suppose that  
$f$ is in $L^\infty(G)$ and
\begin{equation}
\int_G f(x)g(x)dx=0
\label{6.1}
\end{equation}
for all $g$ in $L^1(G)$ such that $\widehat{g}= 0$ on $-B$.
Then $\sigma[f]\subset B$.
\label{prop5.2}
\end{prop}
{\bf Proof.}\quad  Let $\chi_0$ be any
element in $\Gamma\setminus B$.  We will show that
$\chi_0$ is not in the spectrum of $f$ by constructing a 
function $h$ in $L^1(G)$ with $\widehat{h}(\chi_0)\neq 0$ 
and $h*f= 0$.  Let $U$ be an open neighborhood of  
$\chi_0$ not intersecting $B$, and let $h$ be in $L^1(G)$ such that
$\widehat{h}$ is equal to 1 at $\chi_0$ and to 0 outside $U$.
Direct computations show that the Fourier transform
of the function $g:\ \ t\mapsto h(x-t)$, when evaluated at
$\chi\in\Gamma$, gives $\overline{\chi(x)}\widehat{h}(-\chi)$,
and hence it vanishes on $-B$.  
It follows from (\ref{6.1}) that $h*f= 0$, 
which completes the proof.\\


A certain class of $S$-sets, known as the 
Calder\'on sets, or $C$-sets,
is particularly useful to us.
These are defined as follows.
A subset $E$ of $\Gamma$ is 
called a $C$-set if every $f$ in 
$L^1(G)$ with Fourier transform
vanishing on $E$ can be approximated
in the $L^1$-norm by functions of the
form $h*f$ where $h\in L^1(G)$ and $\widehat{h}$ vanishes on
an open set containing $E$.  \\
$C$-sets enjoy the following properties
(see \cite[(39.39)]{hr2} or \cite[Section 7.5]{rudin}).  
\begin{itemize}
\item
Every $C$-set is an $S$-set.
\item
Every closed subgroup of $\Gamma$ is a $C$-set.
\item
The empty set is a $C$-set.    
\item
If the boundary of a set $A$ is a $C$-set,
then $A$ is a $C$-set.  
\item
Finite unions of $C$-sets are $C$-sets.
\end{itemize}

Since closed subgroups are $C$-sets, we conclude that
$\overline{P} \cap \overline{(-P)}$,
and $C_\alpha$, for all $\alpha$,
are $C$-sets.  
>From the definition of $S_{\alpha_0}$, (\ref{alpha0slice}),
and the fact that $C_{\alpha_0}$ is open and closed, it 
follows that the boundary of $S_{\alpha_0}$ is the 
closed subgroup $\psi_{\alpha_0}^{-1}(0)\cap C_{\alpha_0}$.
Hence $S_{\alpha_0}$ is a $C$-set.
For $\alpha\neq \alpha_0$, the set
$S_\alpha$ is open and closed,
and so it has 
empty boundary, and thus it is a $C$-set.  Likewise
$C_\alpha\setminus D_\alpha$ is a $C$-set for all
$\alpha\neq \alpha_0$.  
>From this we conclude that
arbitrary unions of $S_\alpha$ and $C_\alpha\setminus D_\alpha$
are $C$-sets, because an arbitrary union of
such sets, not including the index $\alpha_0$, is open and closed,
and so it is a $C$-set.  \\
We summarize our findings as follows.
\begin{prop}
Suppose that $P$ is a measurable order on $\Gamma$.  
We have:\\
(i)  $\overline{P}$ and $\overline{(-P)}$ are $C$-sets;\\
(ii) $S_\alpha$ is a $C$-set for all $\alpha$;\\
(iii)  arbitrary unions of $S_\alpha$ and $C_\alpha\setminus D_\alpha$
are $C$-sets.
\label{prop5.1}
\end{prop}

As an immediate application, we
have the following characterizations.
\begin{cor}
Suppose that $f$ is in $L^\infty(G)$,
then \\
(i)\  $\sigma[f] \subset S_\alpha$ if and only if 
$\int_G f(x) g(x) d x =0$ for all 
$g\in L^1(G)$ such that $\widehat{g}=0$ on $-S_\alpha$;\\
(ii)\  $\sigma[f] \subset \Gamma\setminus C_\alpha$ if and only if
$\mu_\alpha *f=0$;\\
(iii)\ $\sigma[f] \subset \overline{P}$ if and only if
$f\in H^\infty(G)$.
\label{cor5.3}
\end{cor}

{\bf Proof.}  Assertions (i) and (iii) are clear from
Propositions \ref{prop5.1} and \ref{prop5.2}.  
To prove 
(ii), use Fubini's Theorem to first establish
that for any $g\in L^1(G)$, and any $\mu\in M(G)$, we have
$$\int_G (\mu*f)(t) g(t)dt=\int_G f(t)(\mu*g)(t)dt.$$
Now suppose that 
$\sigma[f] \subset \Gamma\setminus C_\alpha$, and let $g$ be any function in
$L^1(G)$.  From Propositions \ref{prop5.1} and 
\ref{prop5.2}, we have that
$\int_G f g dt = 0$ for all $g$ with Fourier transform vanishing on
$\Gamma\setminus C_\alpha$, equivalently, for all $g=\mu_\alpha*g$.
Hence, $\int_G f (\mu_\alpha*g) dt = \int_G (\mu_\alpha*f) g dt=0$ for all $g$
 in $L^1(G)$, from which it follows that $\mu_\alpha * f=0$.  
 The converse is proved similarly, and we omit the details.

Aiming for a characterization of weakly analytic measures in terms of their spectra, we present one more result.
\begin{prop}
Let $\mu$ be weakly measurable in $M(\Sigma)$.\\
(i)  Suppose that $B$ is a nonvoid closed
subset of $\Gamma$ and ${\rm spec}_T\mu\subset B$.   
Then $\sigma[t\mapsto T_t\mu(A)]\subset B$ for all 
$A\in\Sigma$.\\
(ii)  Conversely, suppose that $B$ is an $S$-set in $\Gamma$
and that $\sigma[t\mapsto T_t\mu(A)] \subset B$ for all
$A\in\Sigma$, then ${\rm spec}_T\mu \subset B$.
\label{prop5.5}
\end{prop}
{\bf Proof.}  We clearly have 
$\cI (\mu)\subset \cI ([t\mapsto T_t\mu(A)])$. 
Hence, ${\rm spec}_T \mu=Z(\cI (\mu))\supset 
Z(\cI([t\mapsto T_t\mu(A)]))=\sigma[t\mapsto T_t\mu(A)],$
and (i) follows.  \\
Now suppose that $B$ is an $S$-set and let $g\in L^1(G)$ 
be such that $\widehat{g}= 0$ on $-B$.  Then, for all
$A\in\Sigma$,
we have from Proposition \ref{prop5.2} (ii):
$$\int_G g (t) T_t\mu (A)d t = 0.$$ 
Equivalently, we
have that
$$\int_G g(-t) T_{-t}\mu(A) dt =0.$$
Since the Fourier transform of the function
$t\mapsto g(-t)$ vanishes on $B$,
we see that
$\cI (\mu) \supset \{f:\ \ \widehat{f}=0\ {\rm on}\ B\}$.
Thus $Z(\cI (\mu))\subset 
Z(\{f:\ \ \widehat{f}=0\ {\rm on} \ B\})=B$,
which completes the proof.

Straightforward applications of Propositions
\ref{prop5.1} 
and \ref{prop5.5} yield the 
desired characterization of weakly analytic measures.
\begin{cor}
Suppose that $\mu\in \cMT$.  Then,\\
(i)  $\mu$ is weakly $T-$analytic if and only if
${\rm spec}_T\mu\subset \overline{P}$ if and only if
$\sigma [t\mapsto T_t\mu (A)]\subset \overline{P}$, for every
$A\in \Sigma$;  \\ 
(ii) $ {\rm spec}_T\mu\subset S_\alpha $ if and only if 
$\sigma[t\mapsto T_t\mu(A)]\subset S_\alpha$
for every $A\in \Sigma$.  \\
(iii) $ {\rm spec}_T\mu\subset C_\alpha $ if and only if 
$\sigma[t\mapsto T_t\mu(A)]\subset C_\alpha$
for every $A\in \Sigma$.
\\
(iv)  $ {\rm spec}_T\mu\subset \Gamma\setminus C_\alpha $
 if and only if 
$\sigma[t\mapsto T_t\mu(A)]\subset \Gamma\setminus  C_\alpha$
for every $A\in \Sigma$. 
\label{cor5.7}
\end{cor}
The remaining results of this section are 
simple properties of measures in $\cMT$ that will be needed later.  
Although the statements are 
direct analogues of classical facts about measures on groups, these
generalization require in some places the sup path
attaining property of $T$.

\begin{prop}
Suppose that $\mu\in \cMT$ and $\nu\in M(G)$.
Then ${\rm spec}_T \nu*_T\mu$ is contained in the support
of $\widehat{\nu}$, and 
${\rm spec}_T \nu*_T\mu\subset {\rm spec}_T \mu$.
\label{proposition5.8}
\end{prop}
{\bf Proof.}  
Given $\chi_0$ not in the support of $\widehat{\nu}$, to
conclude that it is also not in the spectrum of 
$\nu*_T\mu$ it is enough to find a function $f$
in $L^1(G)$ with $\widehat{f}(\chi_0)=1$ and
$f*_T(\nu*_T\mu)=0$.  Simply choose $f$ with Fourier transform
vanishing on the support of $\widehat{\nu}$ and taking value 1
at $\chi_0$.  By Fourier inversion, we have
$f*\nu=0$, and since
$f*_T(\nu*_T\mu)=(f*\nu)*_T\mu$, the 
first part of the proposition follows.
For the second part,
we have $\cI (\mu)\subset \cI (\nu*_T\mu)$,
which implies the desired inclusion.\\

We next prove a property of $L^\infty(G)$ functions
similar to the characterization of $L^1$ functions which are
constant on cosets of a subgroup \cite[Theorem (28.55)]{hr2}.
\begin{prop}\label{proposition5.9}
Suppose that $f$ is in $L^\infty(G)$ and that $\Lambda$ is
an open subgroup of $\Gamma$.  Let $\lambda_0$ denote
the normalized Haar measure on the compact group
$A(G,\Lambda)$, the annihilator in $G$ of
$\Lambda$ (see \cite[(23.23)]{hr1}.  Then, $\sigma[f]\subset \Lambda$
if and only if $f=f*\lambda_0$ a.\ e.\   This is also
the case if and only if $f$ is constant on
cosets of $A(G,\Lambda)$.
\end{prop}
{\bf Proof.}  Suppose that the spectrum of $f$ is contained 
in $\Lambda$.  Since $\Lambda$
is an $S$-set, it follows that $f$
is the weak-* limit of trigonometric polynomials with
spectra contained in $\Lambda$.  Let $\{f_\alpha\}$ be a net
of such trigonometric polynomials converging 
to $f$ weak-*.  Note that, for any $\alpha$, 
we have $\lambda_0*f_\alpha=f_\alpha$.  For $g$ in $L^1(G)$, we have 
$$\lim_\alpha \int_G f_\alpha \overline{g}dx= \int_G f\overline{g}dx.$$
In particular, we have
$$\lim_\alpha \int_G f_\alpha 
(\lambda_0*\overline{g})dx= \int_G f
(\lambda_0*\overline{g})dx,$$
and so
$$\lim_\alpha \int_G (f_\alpha *\lambda_0) 
\overline{g}dx= \int_G (f*\lambda_0)
\overline{g}dx.$$
Since this holds for any $g$ in $L^1(G)$, we conclude that
$\lambda_0 *f_\alpha$ converges weak-* to $\lambda_0*f$.
But $\lambda_0*f_\alpha=f_\alpha$, and $f_\alpha$
converges weak-* to $f$, hence 
$f*\lambda_0=f$.  The remaining assertions of the 
lemma are easy to prove.  We omit the details.

In what follows, we use the symbol $\mu_\alpha$ to denote the
normalized Haar measure on the compact subgroup
$A(G,C_\alpha)$, the annihilator in $G$ of $C_\alpha$.  
This measure is also characterized by its Fourier
transform:
$$\widehat{\mu_\alpha}=1_{C_\alpha}$$
(see \cite[(23.19)]{hr1}).

\begin{cor}
Suppose that $\mu\in\cMT$.  Then,\\
(i)\ ${\rm spec}_T \mu\subset C_\alpha$ if and only if
$\mu=\mu_\alpha *_T \mu$;\\
(ii)\      ${\rm spec}_T \mu\subset \Gamma\setminus C_\alpha$ 
if and only if
$\mu_\alpha *_T \mu=0$.
\label{corollary5.10}
\end{cor}
{\bf Proof.}  (i)  If $\mu=\mu_\alpha*_T\mu$,
then, by Proposition \ref{proposition5.9}, 
$\sigma[t\mapsto \mu_\alpha*_T T_t\mu(A)]\subset C_\alpha$.
Hence by Corollary \ref{cor5.7}, ${\rm spec}_T \mu\subset C_\alpha$.  
For the other direction, suppose that 
${\rm spec}_T \mu\subset C_\alpha$.  Then by Corollary \ref{cor5.7}
we have that the spectrum of the function
$t\mapsto T_t\mu (A)$ is contained in $C_\alpha$
for every $A\in \Sigma$.  By Proposition \ref{proposition5.9}, 
we have that
$$T_t\mu(A)=\int_{G_\alpha}T_{t-y}\mu (A) d\mu_\alpha=
T_t(\mu_\alpha*\mu)(A)$$
for almost all $t\in G$.  Since this holds for all $A\in\Sigma$,
the desired conclusion
follows from Proposition \ref{prop hypa}.\\
Part (ii) follows from Corollary \ref{cor5.3} (ii), Proposition
\ref{prop5.5}(ii), and the fact that $\Gamma\setminus C_\alpha$
is an $S$-set.
\begin{cor}
Suppose that $\mu\in \cMT$ and 
${\rm spec}_T\mu \subset C_\alpha$,
and let $y\in G_\alpha=A(G,C_\alpha)$.  Then
$T_y\mu=\mu$.
\label{cor5.11}
\end{cor}
{\bf Proof.}  For any $A\in \Sigma$, we have from Corollary 
\ref{corollary5.10}
\begin{eqnarray*}
T_y\mu (A)
		&=&
T_y(\mu_\alpha *\mu)(A)=\mu_\alpha*T_y\mu(A)\\
		&=&
\int_{G_\alpha} T_{y-x} \mu (A) d\mu_\alpha (y)\\
		&=&
\int_{G_\alpha}T_{-x}\mu(A) d\mu_\alpha (y) 
=\mu_\alpha *\mu(A)=\mu(A).
\end{eqnarray*}


\section{Homomorphism theorems}

We continue with the notation of the previous section:
$G$ is a locally compact abelian group,
$\Gamma$ the dual group of $G$, $P$ is a measurable
order on $\Gamma$, $T$ is a sup path 
attaining representation of $G$ acting on $M(\Sigma)$.
Associated with $P$ is a collection of homomorphisms
$\psi_\alpha$, as described by Theorem \ref{structureorder}.  Let $\phi_\alpha$ denote the 
adjoint of $\psi_\alpha$.  Thus, $\phi_\alpha$ is a continuous
homomorphism of $\R$ into $G$.
By composing the representation $T$ with the $\phi_\alpha$, we define a new representation $T_{\phi_\alpha}$ of $\R$ acting on $M(\Sigma)$ by:
$t\in \R \mapsto T_{\phi_\alpha(t)}$.  If $\mu$ in $M(\Sigma)$ is weakly measurable with respect to $T$
then $\mu$ is also weakly measurable with respect
to $T_{\phi_\alpha}$.  We will further suppose that
$T_{\phi_\alpha}$ is sup path attaining for each $\alpha$.  This is the case with the representations
of Example \ref{s-set}.

Our goal in this section is to relate the notion of
analyticity with respect to $T$ to the
notion of analyticity with respect to $T_{\phi_\alpha}$.
More generally, 
suppose that $G_1$ and $G_2$ are two
locally compact abelian groups with dual groups $\Gamma_1$ and
$\Gamma_2$, respectively.  Let 
$$\psi: \Gamma_1\rightarrow \Gamma_2$$
be a continuous homomorphism, and
let $\phi:G_2\rightarrow G_1$ denote its adjoint homomorphism.
Suppose
$\nu$ is in $M(G_2)$.  We define a Borel
measure $\Phi(\nu)$ in $M(G_1)$ on the 
Borel subsets $A$ of $G_1$ by:
\begin{equation}
\Phi(\nu)(A)=\int_{G_2} 1_A\circ \phi(t)\,d\nu(t)=
\int_{G_1} 1_A d\Phi(\nu),
\label{continuous-image1}
\end{equation}
where $1_A$ is the indicator function of $A$.
We have
$\|\Phi(g)\|_{M(G_1)}=\|\nu\|_{M(G_2)}$ and, for every Borel measurable
bounded function $f$ on $G_1$, we have
\begin{equation}
\int_{G_1} f d\Phi(\nu)=
\int_{G_2} f\circ \phi(t)\,d \nu (t).
\label{continuous-image2}
\end{equation}
In particular, if $f=\chi$, a character in
$\Gamma$, then
\begin{equation}
\widehat{\Phi(\nu)}(\chi)=\int_{G_1} \overline{\chi} d\Phi(\nu)=\int_{G_2} \overline{\chi}\circ \phi(t)\,d\nu(t)
=\int_{G_2} 
\overline{\psi(\chi)}(t)\, d\nu(t)=\widehat{\nu}(\psi(\chi)),
\label{continuous-image3}
\end{equation}
where $\psi$ is the adjoint homomorphism of $\phi$.
So, 
\begin{equation}
\widehat{\Phi(\nu)}=\widehat{\nu}\circ \psi.
\label{continuous-image4}
\end{equation}

Our first result is a very useful fact
from spectral synthesis of bounded functions.
The proof uses in a crucial way the fact that the representation is sup path attaining, or, more precisely, satisfies the property in Proposition \ref{prop hypa}.

\begin{lemma}
Suppose that $T$ is a sup path attaining representation of $G_1$ acting on $M(\Sigma)$, $\phi$ is a continuous homomorphism of $G_2$ into $G_1$ 
such that $T_\phi$ is a sup path 
attaining representation
of $G_2$.  Suppose that $B$ is a nonempty closed $S$-subset of
$\Gamma_1$ and that $\mu$ is in $M(\Sigma)$ with
$\spec_T\mu\subset B$.  Suppose that
$C$ is an $S$-subset of $\Gamma_2$ and
$\psi(B)\subset C$.
Then $\spec_{T_\phi}\mu\subset C$.
\label{lem3.1}
\end{lemma}

{\bf Proof.}\quad Since $C$ is an $S$-subset of $\Gamma_2$, it is enough to show that for every 
`$A\in \Sigma$, $\spec_{T_\phi}(x\mapsto T_{\phi(x)}\mu(A))\subset C$, by 
Proposition \ref{prop5.5}.
For this purpose, it is enough by 
 \cite[Theorem (40.8)]{hr2}, to show that 
$$g*T_{\phi(\cdot)}\mu(A)=0$$
for every $g$ in $L^1(G_2)$ with
$\widehat{g}=0$ on $C$.  
For $r\in G_2$ and $x\in G_1$, consider the measure
$$T_x(g*_{T_\phi}T_{\phi(r)}\mu)=
g*_{T_\phi}T_{x+\phi(r)}\mu.$$
For $A\in \Sigma$, we have
\begin{eqnarray*}
g*_{T_\phi} T_{x+\phi(r)}\mu(A)&=&
\int_GT_{-t+x}(T_{\phi(r)}\mu)(A)\,d\Phi(g)(t)\\
&=&\Phi(g)*[t\mapsto T_t(T_{\phi(r)}\mu)(A)](x)\\
&=&0
\end{eqnarray*}
for almost all $x\in G_1$.  
To justify the last equality,
we appeal to Proposition \ref{annihilate} and
note that
$\widehat{\Phi(g)}=\widehat{g}\circ \psi$ and so $\widehat{\Phi(g)}=0$ on $B\subset \psi^{-1}(C)$.
Moreover, $\sigma[t\mapsto T_t(T_{\phi(r)}\mu)(A)]
\subset \spec_T(\mu)\subset B$.
Now, using Proposition \ref{prop hypa} and the fact that,
for every $A\in \Sigma$,
$$T_x[g*_{T_\phi} 
T_{\phi(r)}]\mu(A) =g*_{T_\phi}T_{x+\phi(r)}\mu(A)=0,$$
for almost all $x\in G_1$, we conclude that the 
measure $g*_{T_\phi} 
T_{\phi(r)}\mu$ is the zero measure, which completes the proof.

Given ${\cal C}$, a collection of elements in
$L^1(G_1)$ or $M(G_1)$, let 
$$Z({\cal C})=\bigcap_{\delta\in {\cal C}}
\left\{
\chi:\ \widehat{\delta}(\chi)=0\right\}.$$
This is the same notation for the zero set of an ideal
in $L^1(G)$ that we introduced in Section 1.  
Given a set of measures ${\cal S}$ in $M(G_2)$,
let 
$$\Phi({\cal S})=\left\{
\Phi(\nu):\ \nu \in {\cal S}\right\}\subset M(G_1).$$

\begin{lemma}
In the above notation, if 
 $\mu\in M(\Sigma)$ is weakly measurable,
then
$$Z\left( \Phi({\cal I}_{T_\phi}\mu)\right)=
\psi^{-1}\left(    
Z ({\cal I}_{T_\phi}\mu)\right)=
\psi^{-1}\left(    
\spec_{T_\phi}\mu\right).
$$
\label{lem3.2}
\end{lemma}

{\bf Proof.}  It is enough to establish the first equality; the second one follows from definitions.
We have
\begin{eqnarray*}
Z\left( \Phi({\cal I}_{T_\phi}\mu)\right)
&=&
\bigcap_{\delta\in    \Phi({\cal I}_{T_\phi}(\mu))   }
\left\{
\chi\in \Gamma:\ \widehat{\delta}(\chi)=0\right\}\\
&=&
\bigcap_{g \in    {\cal I}_{T_\phi}(\mu)   }
\left\{
\chi\in \Gamma:\ \widehat{\Phi(g)}(\chi)=0\right\}\\
&=&
\bigcap_{g \in    {\cal I}_{T_\phi}(\mu)   }
\left\{
\chi\in \Gamma:\ 
\widehat{g }(\psi(\chi))=0\right\}\\
&=&
\bigcap_{g \in    {\cal I}_{T_\phi}(\mu)   }
\psi^{-1}\left(  Z(g)  \right)\\
&=&
\psi^{-1}\left(
\bigcap_{g \in    {\cal I}_{T_\phi}(\mu)   }
\left(  Z(g)  \right)\right)\\
&=&
\psi^{-1}\left(
Z \left(   {\cal I}_{T_\phi}(\mu)   
  \right)\right)
=\psi^{-1}\left(\spec_{T_\phi}\mu\right).
\end{eqnarray*}

\begin{lemma}
Suppose that $C$ is a nonempty closed $S$-subset of
$\Gamma_2$ and that $\psi^{-1}(C)$ is an 
$S$-subset of $\Gamma_1$.  Suppose that $\mu$ is in $M(\Sigma)$ and 
$\spec_{T_\phi}(\mu)\subset C$.  
Then $\spec_T\mu\subset \psi^{-1}(C)$.
\label{lem3.3}
\end{lemma}

{\bf Proof.}\quad 
We will use the notation of Lemma \ref{lem3.2}.
If $f\in \cI_{T_\phi}(\mu)$ and $t\in G_1$, 
then $f\in \cI_{T_\phi}(T_t \mu)$.  So,
 for $A\in\Sigma$, we have
$f*_{T_\phi}(T_t\mu)(A)=0$.  
But
\begin{eqnarray*}
f*_{T_\phi}(T_t\mu)(A)&=&
\int_\R T_{t-\phi(x)}\mu(A)f(x)\, dx\\
&=&
\int_G T_{t- x }\mu(A)\, d\Phi(f),
\end{eqnarray*}
where $\Phi(f)$ is the homomorphic image
of the measure $f(x)\,dx$.
Hence,
$\Phi(f) \in \cI_T ([t\mapsto T_t\mu(A)])$,
and so $\Phi(\cI_{T_\phi}(\mu))\subset 
\cI_T([t\mapsto T_t\mu(A)])$, which implies that
  $$Z\left(\Phi(\cI_{T_\phi}(\mu))\right)\supset 
Z\left(\cI_T([t\mapsto T_t\mu(A)])\right)
=\spec_T(t\mapsto T_t\mu(A)).$$
By Lemma \ref{lem3.2}, 
$$Z\left(\Phi(\cI_{T_\phi}(\mu))\right)=
\psi^{-1}\left(
\spec_{T_\phi}\mu
\right)\subset \psi^{-1}(C).$$
Hence, $\spec_T(t\mapsto T_t\mu(A))\subset \psi^{-1}(C)$
for all $A\in\Sigma$, which by Proposition \ref{prop5.5}  implies that
$\spec_T(\mu)\subset  \psi^{-1}(C)$.

Taking $G_1=G,\ G_2=\R$ and $\psi=\psi_\alpha$ to be one of the 
homomorphisms in Theorem \ref{structureorder}, 
and using the fact that $[0,\infty[$, $S_\alpha$,
$C_\alpha\setminus D_\alpha$ are all $S$-sets,
we obtain  useful relationships between different 
types of analyticity.

\begin{thm}\label{equiv-def}
Let $G$ be a locally compact abelian group with ordered dual group
$\Gamma$, and let $P$ denote a measurable order on $\Gamma$.  
Suppose that $T$ is a sup path attaining representation
of $G$ by isomorphisms of $M(\Sigma)$, such that $T_{\phi_\alpha}$ is sup path attaining, where 
$\phi_\alpha$ is as in Theorem \ref{structureorder}.\\
(i)  If $\mu\in M(\Sigma)$ and 
$\spec_T(\mu)\subset C_\alpha\setminus D_\alpha$.
Then
$$\spec_T(\mu)\subset S_\alpha\Leftrightarrow
\spec_{T_{\phi_\alpha}}(\mu)\subset [0,\infty[.$$
(ii)  
If $\mu\in M(\Sigma)$ and 
$\spec_T(\mu)\subset C_{\alpha_0}$.
Then
$$\spec_T(\mu)\subset S_{\alpha_0}\Leftrightarrow
\spec_{T_{\phi_{\alpha_0}}}(\mu)\subset [0,\infty[.$$
\end{thm}

We can use the representation $T_\phi$
to convolve a measure $\nu\in M(G_2)$
with $\mu\in M(G_1)$:
$$\nu*_{T_\phi}\mu(A)=\int_{G_2}T_{-\phi(x)}\mu(A)d\nu(x)=\int_{G_2}\mu(A-\phi(x))\,d\nu(x) ,$$
for all Borel $A\subset G_1$.  

Alternatively, we can convolve
$\Phi(\nu)$ 
 in the usual sense of \cite[Definition 19.8]{hr1} 
with $\mu$ to yield
another measure in $M(G_1)$, defined on the Borel subsets of $G_1$ by
$$\Phi(\nu)*\mu(A)=\int_{G_1}\int_{G_1}1_A(x+y)d\Phi(\nu)(x)d\mu(y).$$
Using (\ref{continuous-image2}), we find that
\begin{eqnarray*}
\Phi(\nu)*\mu(A)        &=&
\int_{G_1}\int_{G_2}1_A(\phi(t)+y)d\nu(t)d\mu(y)\\
			&=&
\int_{G_2}\mu(A-\phi(t))d\nu(t)
=\nu *_{T_\phi}\mu(A).
\end{eqnarray*}
Thus, 
\begin{equation}
\Phi(\nu)*\mu=\nu*_{T_\phi}\mu.
\label{7.feb.95.2}
\end{equation}
We end the section 
with homomorphism theorems, which complement the well-known
homomorphism theorems for $L^p$-multipliers 
(see Edwards and Gaudry \cite[Appendix B]{eg}).
In these theorems, we let $G_1$ act on $M(G_1)$ by translation.  That is, if $\mu\in M(G_1)$, $x\in G_1$,
and $A$ is a Borel subset of $G_1$, then
$$T_x\mu(A)=\mu(A+x).$$
Let $\phi:\ G_2\rightarrow G_1$ be a continuous homomorphism.
By Example \ref{exhypa}, $T$ and $T_\phi$ are sup path attaining.
(Recall that if $t\in G_2$, $\mu\in M(G_1)$,
then $T_{\phi(t)}\mu(A)=\mu(A+\phi(t))$.)
A simple exercise with definitions shows that
for $\mu\in M(G_1)$ 
$$\spec_T\mu=\supp \widehat{\mu}.$$

\begin{thm}
Suppose that $\Gamma_1$ and $\Gamma_2$ contain 
measurable orders $P_1$ and $P_2$, respectively, and
$\psi:\ \Gamma_1\rightarrow \Gamma_2$ is a continuous,
order-preserving homomorphism (that is, $\psi(\overline{P_1})\subset
\overline{P_2}$).  
Suppose that there is a positive 
constant $N(\nu)$ such that for all $f\in H^1(G_2)$
\begin{equation}
\|\nu*f\|_1\leq N(\nu)\|f\|_1.
\label{hom1eq}
\end{equation}
Then
\begin{equation}
\|\Phi(\nu)*\mu\|\leq N(\nu)\|\mu\|
\label{hom2eq1}
\end{equation}
for all Borel measures in $M(G_1)$ such that
$\widehat{\mu}$ is supported in $\overline{P_1}$.
\label{homth2}
\end{thm}

{\bf Proof.}\quad We have $\Phi(\nu)*\mu=\nu*_{T_\phi}\mu$.
Also $\overline{P_2}$ is a ${\cal T}$-set.
So (\ref{hom2eq1}) will follow from Theorem 
\ref{trans-thm}
once we show that $\spec_{T_\phi}\mu\subset \overline{P_2}$.  For that purpose, we use Lemma \ref{lem3.1}.
We have 
$$\spec_T\mu =\supp \widehat{\mu}\subset \overline{P_1},$$
and $\psi(\overline{P_1})\subset \overline{P_2}$ is an $S$-set.  Hence 
$\spec_{T_\phi}\mu\subset \overline{P_2}$
by Lemma \ref{lem3.1}.

The following special case of Theorem \ref{homth2} deserves a separate 
statement.

\begin{thm}
With the above notation, suppose that there is a positive 
constant $N(\nu)$ such that for all $f\in H^1(G_2)$
\begin{equation}
\|\nu*f\|_1\leq N(\nu)\|f\|_1.
\label{hom1eq1}
\end{equation}
Then for all $f\in H^1(G_1)$ we have
\begin{equation}
\|\Phi(\nu)*f\|_1\leq N(\nu)\|f\|_1.
\label{hom1eq2}
\end{equation}
\label{homth1}
\end{thm}

\begin{thm}
Suppose that there is a positive 
constant $N(\nu)$ such that for all $f\in H^1(\R)$
\begin{equation}
\|\nu*f\|_1\leq N(\nu)\|f\|_1.
\label{hom2eq11}
\end{equation}
Then for all $\mu\in M(G)$
with support of $\widehat{\mu}$ contained in
$C_\alpha\setminus D_\alpha$, where $\alpha<\alpha_0$, 
we have
\begin{equation}
\|\Phi_\alpha(\nu)*\mu\|_1\leq N(\nu)\|\mu\|.
\label{hom2eq2}
\end{equation}
\label{varianthomth2}
\end{thm}

{\bf Proof.}\quad The proof is very much like the 
proof of Theorem \ref{homth2}.  
We have $\Phi_\alpha (\nu)*\mu=\nu*_{T_\phi}\mu$.
Apply Theorem \ref{trans-thm}, taking into consideration that
$$\spec_T\mu=\supp \widehat{\mu}\subset C_\alpha\setminus D_\alpha$$
is an $S$-set and so 
$$\spec_{T_{\phi_\alpha}}\mu\subset \psi_\alpha(C_\alpha\setminus D_\alpha)\subset [0,\infty[.$$


\section{Decomposition of Analytic Measures}

Define measures
$\mu_{\alpha_0}$ and $d_\alpha$ by their Fourier transforms: 
$\widehat{\mu_{\alpha_0}}=1_{C_{\alpha_0}}$,
and $\widehat{d_\alpha}=1_{C_\alpha\setminus D_\alpha}$.
Then we have the following decomposition theorem.

\begin{thm}
\label{decomp of measures}
Let $G$ be a locally compact abelian group with an ordered dual
group $\Gamma$.
Suppose that
$T$ is a sup path attaining representation of $G$ in $M(\Sigma)$.
Then for any weakly analytic measure $\mu \in M(\Sigma)$ we have
that the set of $\alpha$ for which $d_\alpha *_T\mu \ne 0$ is countable,
and that
\begin{equation}
\mu = \mu_{\alpha_0}*_T\mu+ \sum_\alpha d_\alpha *_T\mu ,
\end{equation}
where the right side converges unconditionally in norm
in $M(\Sigma)$.
Furthermore, there is a positive constant $c$, depending only upon
$T$, such that for any signs $\epsilon_\alpha = \pm 1$ we have
\begin{equation}
\left\| \sum_\alpha \epsilon_\alpha d_\alpha*_T\mu \right\| \le c \|\mu\| .
\end{equation}
\end{thm}
One should compare this theorem
 to the well-known results
from Littlewood-Paley theory on $L^p(G)$, where $1<p<\infty$ (see Edwards and Gaudry \cite{eg}).
For $L^p(G)$ with 
$1<p<\infty$, 
it is well-known that the subgroups $(C_\alpha)$
form a Littlewood-Paley decomposition of the group $\Gamma$, which means that the martingale difference series 
$$f= \mu_{\alpha_0}*f+ \sum_\alpha d_\alpha *f $$
converges unconditionally in $L^p(G)$ to $f$.
Thus, Theorem \ref{decomp of measures} above may be considered as an extension
of Littlewood-Paley Theory to spaces of analytic measures.

The next result, crucial to our proof of Theorem~\ref{decomp of measures},
is already known in the case that $G = \T^n$ with the lexicographic order
on the dual.  This is due to Garling \cite{gar}, and is a modification
of the celebrated inequalities of Burkholder.
Our result
can be obtained directly from the result in \cite{gar} 
by combining the techniques of
\cite{ams3} with the homomorphism theorem~\ref{homth2}.  However, we shall
take a different approach, in effect reproducing Garling's proof in this
more general setting.

\begin{thm}
\label{ucc of h1 fts}
Suppose that $G$ is a locally compact group with ordered dual $\Gamma$.  
Then for
$f\in H^1(G)$, for any set $\{\alpha_j\}_{j=1}^n$ of indices
less than $\alpha_0$, and
for any numbers $\epsilon_j \in \{0,\pm 1\}$ ($1 \le j \le n$),
there is an absolute constant $a>0$ such that
\begin{equation}\label{34}
\left\|
\sum_{j=1}^n
\epsilon_j d_{\alpha_j} *f
\right\|_1
\leq a \|f\|_1.
\end{equation}
Furthermore, 
\begin{equation}
\label{ucc of h1 equation}
f = \mu_{\alpha_0}*f + \sum_\alpha d_\alpha*f,
\end{equation}
where the right hand side converges unconditionally in the norm topology on
$H^1(G)$.
\end{thm}
{\bf Proof.}
The second part of Theorem~\ref{ucc of h1 fts} follows easily
from the first part and Fourier inversion.

Now let us show that if we have the result for compact $G$, then
we have it for locally compact $G$.
Let $\pi_{\alpha_0} :\ \Gamma\rightarrow 
\Gamma/C_{\alpha_0}$ denote the quotient homomorphism
of $\Gamma$ onto the discrete group $\Gamma/C_{\alpha_0}$ 
(recall that $C_{\alpha_0}$ is open),
and define a measurable order on 
$\Gamma/C_{\alpha_0}$ to be $\pi_{\alpha_0}(P)$.  
By Remarks \ref{remarkstructureorder} (c),
the decomposition of the group $\Gamma/C_{\alpha_0}$ that we get
by applying Theorem \ref{structureorder}
to that group, is precisely the 
image by $\pi_{\alpha_0}$ of the decomposition of the group $\Gamma$.  
Let $G_0$ denote the 
compact dual group of $\Gamma/C_{\alpha_0}$.  Thus if
Theorem \ref{ucc of h1 fts} holds for $ H^1(G_0)$, then
applying Theorem \ref{homth2}, we see that 
Theorem \ref{ucc of h1 fts} holds for $G$.

Henceforth, let us suppose that $G$ is compact.
We will suppose that the Haar measure on $G$ is normalized, so
that $G$ with Haar measure is a probability space.  

Since each one of the subgroups $C_\alpha$, and
$D_\alpha$ ($\alpha<\alpha_0$) is open, it follows that 
their annihilators in $G$,
$G_\alpha=A(G,C_\alpha)$, and $A(G,D_\alpha)$, are compact.
Let $\mu_\alpha$ and $\nu_\alpha$ denote 
the normalized Haar measures on $A(G,C_\alpha)$
and $A(G,D_\alpha)$, respectively.
We have 
$\widehat{\mu}_\alpha=1_{C_\alpha}$ (for all
$\alpha$),
and 
$\widehat{\nu}_\alpha=1_{D_\alpha}$ 
(for all $\alpha\neq \alpha_0$), so that
$d_\alpha=\mu_\alpha -\nu_\alpha$.

For each $\alpha$,
let ${\cal B}_\alpha$ denote the 
$\sigma$-algebra of subsets of $G$ of the form
$A+G_\alpha$, where $A$ is a Borel subset of $G$.  
We have ${\cal B}_{\alpha_1}\subset
{\cal B}_{\alpha_2}$, whenever $\alpha_1>\alpha_2$. 
It is a simple matter to see that for 
$f\in L^1(G)$, the conditional expectation of
$f$ with respect to ${\cal B}_\alpha$
is equal to
$\mu_\alpha*f$ 
(see \cite[Chapter 5, Section 2]{eg}).

We may suppose without
loss of generality that
$\alpha_1>\alpha_2>\ldots>\alpha_n$.
Thus the $\sigma$-algebras ${\cal B}_{\alpha_k}$ form a filtration, and the
sequence $(d_{\alpha_1}*f, d_{\alpha_2}*f,\ldots,d_{\alpha_n}*f)$ 
is a martingale
difference sequence with respect to this filtration.

In that case, we have the following result
of Burkholder \cite[Inequality (1.7)]{bur}, and \cite{bur1}.  If
$0<p<\infty$, then there is a positive constant $c$, depending only
upon $p$, such that
\begin{equation}
\label{burkholder's inequality}
\left\| \sup_{1 \le k \le n} 
\left( \sum_{j=1}^k \epsilon_j d_{\alpha_j}*f\right) \right\|_p
\le
c
\left\| \sup_{1 \le k \le n} 
\left( \sum_{j=1}^k d_{\alpha_j}*f\right) \right\|_p.
\end{equation}

\begin{lemma}
For any index $\alpha$, $0<p<\infty$, and $f\in H^1(G)\cap L^p(G)$, we have  
almost everywhere on $G$
\begin{equation}
\left|
\mu_\alpha *f
\right|^p
\leq
\mu_\alpha*\left| f\right|^p,
\label{eq1.95}
\end{equation}
where $\mu_\alpha$ is the normalized Haar measure on the compact subgroup
$G_\alpha=A(G,C_\alpha)$.
\label{improved jensen}
\end{lemma}
{\bf Proof.}
The dual group of $G_\alpha$ is
$\Gamma/C_\alpha$ and can be ordered by the set
$\pi_\alpha (P)$, where 
$\pi_\alpha$ is the natural homomorphism of $\Gamma$ onto 
$\Gamma/C_\alpha$.

Next, by convolving 
with an approximate identity
for $L^1(G)$
consisting of trigonometric
polynomials, we may assume  
that $f$ is a 
trigonometric polynomial.
Then we see that for each $x \in G$ that 
the function $y\mapsto f(x+y)$, $y\in G_\alpha$, is 
in $H^1(G_\alpha)$.  To verify this, it is sufficient
to consider the case when
$f$ is a character in $P$.  Then
$$ f(x+y)= f(x) \pi_\alpha (f)(y), $$
and by definition $\pi_\alpha(f)$ is in $H^1(G_\alpha)$.

Now we have the following generalization of Jensen's Inequality,
due to Helson and Lowdenslager 
\cite[Theorem 2]{hl1}.  An independent proof based on the ideas of this 
section is given in \cite{ams3}.  For all $g\in H^1(G)$
\begin{equation}
\left|\int_G g(x) d x\right| \leq \exp \int_G\log
|g(x)|d x.
\label{jensen's inequality}
\end{equation}
Apply (\ref{jensen's inequality}) to $y\mapsto f(x+y)$, $y\in G_\alpha$
to obtain
$$
\left|\int_{G_\alpha} f(x+y) d \mu_\alpha(y)\right| 
\leq \exp \int_{G_\alpha}\log
|f(x+y)|d \mu_\alpha(y).
$$
Extending the integrals to the whole of $G$ (since $\mu_\alpha$ is 
supported on $G_\alpha$),
raising both sides to the 
$p$th power, and then applying
the usual Jensen's inequality for the 
logarithmic function on finite measure spaces, 
we obtain
\begin{eqnarray*}
\left|\int_G f(x+y) d \mu_\alpha(y)\right|^p 
	&\leq& 
\exp \int_G \log
|f(x+y)|^p d \mu_\alpha(y)\\
	&\leq&
\int_G |f(x+y)|^pd \mu_\alpha(y).
\end{eqnarray*}
Changing $y$ to $-y$, we obtain the desired inequality.\\

Let us continue with the proof of Theorem~\ref{ucc of h1 fts}.
We may suppose that $f$ is a mean zero
trigonometric polynomial, and that
the spectrum of $f$ is contained in 
$\bigcup_{j=1}^n C_{\alpha_j} \setminus D_{\alpha_j} $, that is to say
$$ f = \sum_{j=1}^n d_{\alpha_j} * f .$$
By Lemma~\ref{improved jensen}, we have that
\begin{eqnarray}
\sup_{1\leq k\leq n} 
\left|
\mu_{\alpha_k}*f
\right|   
		&=&
\left(
\sup_{1\leq k\leq n} 
\left|
\mu_{\alpha_k}*f
\right|^{1/2}\right)^2   \nonumber\\
		&\leq&
\left(
\sup_{1\leq k\leq n} 
\mu_{\alpha_k}*| f|^{1/2}
\right)^2.   
\label{ucc 2}
\end{eqnarray}
Also, we have that 
$(\mu_{\alpha_j}*|f|^{1/2})_{j=1}^n$
is a martingale with respect to the filtration $({\cal B}_j)_{j=1}^n$.
Hence, by Doob's Maximal
Inequality \cite[Theorem (3.1), p. 317]{doob} we have that
\begin{eqnarray}
\left\| 
\sup_{1\leq k\leq n'} 
\mu_{\beta_k}*| f|^{1/2}
\right\|_2^2
			&\leq&
4
\left\|
\mu_{\beta_{n'}}*| f|^{1/2}
\right\|_2^2              \nonumber\\
			&\leq& 
4 \left\| |f|^{1/2}\right\|_2^2 = 4\|f\|_1.
\label{ucc 3}
\end{eqnarray}
The desired inequality follows now upon combining
Burkholder's Inequality (\ref{burkholder's inequality})
with (\ref{ucc 2}), and (\ref{ucc 3}).\\

\noindent
{\bf Proof of Theorem~\ref{decomp of measures}.}  
Transferring inequality (\ref{34}) by using Theorem \ref{trans-thm}, we obtain that 
for any set $\{\alpha_j\}_{j=1}^n$ of indices
less than $\alpha_0$, and
for any numbers $\epsilon_j \in \{0,\pm 1\}$ ($1 \le j \le n$),
there is a positive constant $c$, depending only upon the
representation $T$, such that
\begin{equation}
\left\|
\sum_{j=1}^n
\epsilon_j d_{\alpha_j} *_T \mu
\right\|
\leq c \|\mu\|.
\end{equation}

Now suppose that $\{\alpha_j\}_{j=1}^\infty$ is a countable collection of
indices less than $\alpha_0$.  
Then by Bessaga and Pe\l czy\'nski \cite{bp}, the series
$\sum_{j=1}^\infty d_{\alpha_j} *_T \mu$ is unconditionally convergent.
In particular, for any $\delta>0$, for only finitely many $k$ do we have
that $\| d_{\alpha_k} *_T \mu \| > \delta$.  Since this is true for all
such countable sets, we deduce that the set of $\alpha$ for which
$ d_\alpha *_T \mu \ne 0$ is countable.

Hence we have that 
$\sum_\alpha d_{\alpha} *_T \mu$ is unconditionally convergent
to some measure, say $\nu$.
Clearly $\nu$ is weakly measurable.  To prove that
$\mu=\nu$, it is enough by Proposition \ref{prop hypa} to show that
for every $A\in\Sigma$, we have
$T_t\mu(A)=T_t\nu(A)$
for almost all $t\in G$.

We first note that since for every $f\in L^1(G)$
the series $\mu_{\alpha_0}*f+ \sum_\alpha d_\alpha *f$ converges to $f$
in $L^1(G)$, it follows that, for every $g\in L^\infty(G)$,
the series $\mu_{\alpha_0}*g+ \sum_\alpha d_\alpha *g$
converges to $g$ in the weak-* topology of $L^\infty(G)$.
Now on the one hand, for $t\in G$ and $A\in \Sigma$, we have
$\mu_{\alpha_0}*_TT_t\mu(A)+ \sum_\alpha d_\alpha *_T T_t\mu(A)=T_t\nu(A)$, 
because of the (unconditional) convergence of
the series $\mu_{\alpha_0}*_T\mu+ \sum_\alpha d_\alpha *_T\mu$ to $\nu$.
On the other hand, by considering the $L^\infty(G)$ function
 $t\mapsto T_t(A)$, we have that
$\mu_{\alpha_0}*_TT_t\mu(A)+ \sum_\alpha d_\alpha *_T T_t\mu(A)=
\mu_{\alpha_0}*T_t\mu(A)+ \sum_\alpha d_\alpha * T_t\mu(A)=T_t\mu(A)$,
weak *.  Thus $T_t\mu(A)=T_t\nu(A)$ for almost all $t\in G$, and the
proof is complete.


\section{Generalized F. and M. Riesz Theorems}


Throughout this section, 
we adopt the notation of Section 5, specifically, 
the notation and assumptions of 
Theorem \ref{decomp of measures}.

Suppose that $T$ is a sup path attaining representation of $\R$ by isomorphisms of $M(\Sigma)$.
In \cite{amss}, we proved the following 
result concerning bounded operators $\cP$
from $M(\Sigma)$ into $M(\Sigma)$
that commute with the representation $T$ in the following sense:
$$\cP\circ T_t=T_t\circ \cP$$
for all $t\in \R$.  
\begin{thm}
Suppose that $T$ is a representation of $\R$ that is sup path
attaining,
and that $\cP$ commutes with $T$.
Let $\mu\in M(\Sigma)$ be weakly analytic.
Then $\cP \mu$ is also weakly analytic.
\label{caseofR}
\end{thm}

To describe an interesting application
of this theorem from \cite{amss}, let us recall the following.

\begin{defin}
Let $T$ be a sup path attaining 
representation of $G$ in $M(\Sigma)$.
A weakly measurable $\sigma$ in $M(\Sigma)$ is
called quasi-invariant if 
$T_t\sigma$ and $\sigma$
are mutually absolutely continuous for all $t\in G$.  Hence
if $\sigma$ is quasi-invariant
and $A\in \Sigma$, then 
$|\sigma|(A)=0$ if and only if $|T_t(\sigma)|(A)=0$
for all $t\in G$.
\label{qi}
\end{defin}


Using Theorem \ref{caseofR} we obtained in \cite{amss} the following extension of results of 
de Leeuw-Glicksberg \cite{deleeuwglicksberg} and 
Forelli \cite{forelli}, concerning quasi-invariant measures.

\begin{thm}
Suppose that $T$ is a sup path attaining representation
of $\R$ by isometries of $M(\Sigma)$.  Suppose 
that $\mu\in M(\Sigma)$ is weakly analytic, and
$\sigma$ is quasi-invariant.  Write
$\mu=\mu_a+\mu_s$ for the Lebesgue decomposition of $\mu$
with respect to $\sigma$.  Then both
$\mu_a$ and $\mu_s$ are weakly analytic.  In particular,
the spectra of $\mu_a$ and $\mu_s$ are
contained in $[0,\infty)$.
\label{lebesgue-decomp-forR}
\end{thm}

Our goal in this section is to extend Theorems \ref{caseofR}
above to representations of a locally compact abelian group $G$ with ordered dual group $\Gamma$.
More specifically, we want to prove the following theorems.

\begin{thm}
\label{application1}
Suppose that $T$ is a sup path attaining 
representation of $G$ by isomorphisms of $M(\Sigma)$ such that 
$T_{\phi_\alpha}$ is sup path attaining for each $\alpha$.
Suppose that $\cP$ commutes with $T$ in the sense
that 
$$\cP\circ T_t=T_t\circ \cP$$
for all $t\in G$.  
Let $\mu\in M(\Sigma)$ be weakly analytic.
Then $\cP \mu$ is also weakly analytic.
\end{thm}

As shown in \cite[Theorem (4.10)]{amss} 
for the case $G=\R$, 
an immediate corollary of Theorem \ref{application1}
is the following result.

\begin{thm}
\label{application2}
Suppose that $T$ is a sup path attaining representation
of $G$ by isometries of $M(\Sigma)$, such that 
$T_{\phi_\alpha}$ is sup path attaining for each $\alpha$.
  Suppose 
that $\mu\in M(\Sigma)$ is weakly analytic with respect to $T$, and
$\sigma$ is quasi-invariant with respect to $T$.  Write
$\mu=\mu_a+\mu_s$ for the Lebesgue decomposition of $\mu$
with respect to $\sigma$.  Then both
$\mu_a$ and $\mu_s$ are weakly analytic with respect to $T$.  In particular,
the $T$-spectra of $\mu_a$ and $\mu_s$ are
contained in $\overline{P}$.
\end{thm}

{\bf Proof of Theorem \ref{application1}.}\quad
Write
$$\mu=\mu_{\alpha_0}*_T\mu +\sum_\alpha d_\alpha *_T\mu,$$
as in (\ref{decomp of measures}), where the series converges unconditionally
in $M(\Sigma)$.
Then
\begin{equation}
\label{-3}
\cP\mu=\cP(\mu_{\alpha_0}*_T\mu) +\sum_\alpha 
\cP(d_\alpha *_T\mu).
\end{equation}
It is enough to show that the $T$-spectrum of
each term is contained in $\overline{P}$.
Consider the measure 
$\mu_{\alpha_0}*_T\mu$.  We have 
$\spec_T(\mu_{\alpha_0}*_T\mu)\subset S_{\alpha_0}$.
Hence by Theorem \ref{equiv-def}, 
$\mu_{\alpha_0}*_T\mu$ is $T_{\phi_{\alpha_0}}$-analytic.  Applying Theorem \ref{caseofR}, we see that
\begin{equation}
\label{-2}
\spec_{T_{\phi_{\alpha_0}}} (\cP (\mu_{\alpha_0}*_T\mu))\subset [0,\infty[.
\end{equation}
Since $\cP$ commutes with $T$, it is easy to see  
from Proposition \ref{proposition5.9} and Corollary \ref{corollary5.10} that
$$\spec_T (\cP (\mu_{\alpha_0}*_T\mu))\subset C_{\alpha_0}.$$
Hence by (\ref{-2}) and Theorem \ref{equiv-def}, 
$$\spec_T (\cP (\mu_{\alpha_0}*_T\mu))\subset S_{\alpha_0},$$
which shows the desired result for the first term
of the series in 
(\ref{-3}).  The other terms of the series 
(\ref{-3}) are handled similarly.
\\

\noindent
{\bf Acknowledgments}  The second author is 
grateful for financial support from the National Science Foundation (U.S.A.) and the Research Board of the 
University of Missouri.

\end{document}